\documentclass[12pt]{article}
\usepackage{amssymb}
\usepackage{amsfonts}
\usepackage{times}
\usepackage{mathptmx}
\usepackage{amsmath}
\usepackage[usenames]{color}
\usepackage{mathrsfs}
\usepackage{amsfonts}
\usepackage{amssymb,amsmath}
\usepackage{CJK}
\usepackage{cite}
\usepackage{cases}
\usepackage{amsthm}

\def\cl{\centerline}

\def\a{\alpha}

\def\b{\beta}
\def\vs{\vspace*}

\def\L{\mathscr{W}}

\def\Z{\mathbb{Z}}
\def\N{\mathbb{N}}

\def\C{\mathbb{C}}

\def\QED{\hfill$\Box$}
\def\ni{\noindent}

\numberwithin{equation}{section}
\newtheorem{theo}{Theorem}[section]
\newtheorem{defi}[theo]{Definition}

\newtheorem{lemm}[theo]{Lemma}
\newtheorem{prop}[theo]{Proposition}

\begin{document}
\begin{center}
{\bf\large Loop Virasoro Lie Conformal Algebra}

\end{center}

\cl{Henan Wu\,$^{1,a)}$, Qiufan Chen\,$^2$, Xiaoqing Yue\,$^{3}$}
\footnote {$^{a)}\,$ Author to whom correspondence should be addressed. Electronic mail: wuhenanby@163.com.
}
\cl{\small\it $^{1,2,3}$Department of Mathematics, Tongji University, Shanghai,
200092, China}

\vs{8pt}

{\small
\parskip .005 truein
\baselineskip 3pt \lineskip 3pt
\begin{abstract}
The Lie conformal algebra of loop Virasoro algebra, denoted by
 $\mathscr{CW}$, is introduced in this paper. Explicitly,  $\mathscr{CW}$ is a Lie conformal algebra with $\C[\partial]$-basis $\{L_i\,|\,i\in\Z\}$ and $\lambda$-brackets $[L_i\, {}_\lambda \, L_j]=(-\partial-2\lambda) L_{i+j}$.
Then conformal derivations of  $\mathscr{CW}$ are determined. Finally, rank one conformal modules and $\Z$-graded free intermediate series modules over  $\mathscr{CW}$ are classified.
\end{abstract}


\parskip .001 truein\baselineskip 6pt \lineskip 6pt

\section{INTRODUCTION}
The notion of Lie conformal algebras, introduced in Ref.$7$, encodes an axiomatic description of the operator
product expansion of chiral fields in conformal field theory. Conformal module is a basic tool for the construction of free field realization
of infinite dimensional Lie (super)algebras in conformal field theory.
In recent years, the structure theory, representation theory and cohomology theory of Lie conformal
algebras have been  extensively studied by many scholars.
For example, a finite simple conformal algebra was proved to be isomorphic to either the conformal Virasoro algebra or  the current conformal algebra associated with a simple finite dimensional Lie algebra in Ref.$3$. Finite irreducible conformal modules over the conformal Virasoro algebra were determined in Ref.$2$.
The cohomology theory of conformal algebras was developed in Ref.$1$. The low dimensional cohomologies of the infinite rank general Lie
conformal algebras $gc_N$ with trivial coefficients were computed in Ref.$9$. Filtered Lie conformal
algebras whose associated graded algebras are isomorphic to that of general conformal algebra $gc_1$
were investigated in Ref.$11$. Two new nonsimple conformal algebras associated with the
$Schr\ddot{o}dinger-Virasoro$ Lie algebra and the extended $Schr \ddot{o} dinger-Virasoro$ Lie algebra were constructed in Ref.$10$.
The Lie conformal algebra of a block type was introduced and free intermediate series modules were classified in Ref.$4$.

In this paper, we would like to consider a Lie conformal algebra related to the so-called {\it loop Virasoro algebra}, whose representation theory was studied  in Ref.$5$. The {\it loop Virasoro algebra} $\L$ is defined to be a Lie algebra with basis $\{L_{\a,i}\,|\,\a, i\in\Z\}$ and Lie brackets given by
\begin{equation}\label{1.1}
[L_{\a,i},L_{\b,j}]=(\b-\a)L_{\a+\b,i+j}. 
\end{equation}
The loop Virasoro algebra plays an important role in the structure
theory of infinite dimensional Lie algebras.$^6$
The Lie conformal algebra of $\L$, denoted by $\mathscr{CW}$, is constructed in Section $2$.
As one can see, it is a Lie conformal algebra with $\C[\partial]$-basis $\{L_i\,|\,i\in\Z\}$ and $\lambda$-brackets
\begin{equation}\label{1.3}
[L_i\, {}_\lambda \, L_j]=(-\partial-2\lambda) L_{i+j}.
\end{equation}
We should remark that the conformal subalgebra $\C[\partial]L_0$ is isomorphic to the well known conformal Virasoro algebra.

This paper is organized as follows. In Section $2$, we recall some basic definitions of Lie conformal algebras. In Section $3$, conformal derivations of $\mathscr{CW}$ are determined. Furthermore, rank one conformal modules and $\Z$-graded free intermediate series modules over $\mathscr{CW}$ are classified in Section $4$ and Section $5$,  respectively.

\section{PRELIMINARIES}

We recall some definitions related to Lie conformal algebras in this section.$^{3,7,8}$

A formal distribution (usually called a field by physicists) with coefficients in a complex
vector space $U$ is a series of the following form:
\begin{equation*}
a(z)=\sum_{i\in\Z}a_{(i)} z^{-i-1},
\end{equation*}
where $z$ is an indeterminate and $a_{(i)}\in U$.
Denote by $U[[z,z^{-1}]]$ the space of formal distribution with coefficients in $U$.
The space $U[[z,z^{-1},w,w^{-1}]]$ is defined in a similar way. A formal distribution $a(z,w)\in U[[z,z^{-1},w,w^{-1}]]$ is called {\it local} if $(z-w)^N a(z,w)=0$ for some $N\in\Z^+$.
Let $g$ be a Lie algebra.
Two formal distributions $a(z), b(z)\in g[[z,z^{-1}]]$ are
called {\it pairwise local} if $[a(z),b(w)]$ is local in $g[[z,z^{-1},w,w^{-1}]]$.
\begin{defi}\rm
A family $F$ of pairwise local formal distributions, whose coefficients span $g$, is called
a formal distributions Lie algebra of $g$. In such a case, we say that the family $F$ spans $g$. We will write $(g, F)$ to
emphasize the dependence on $F$.
\end{defi}

Define the {\it formal delta distribution} to be
\begin{equation*}
\delta(z,w)=\sum_{i\in\Z}z^i w^{-i-1}.
\end{equation*}
The following proposition describes an equivalent condition for a formal
distribution to be local.$^7$
\begin{prop}\label{local}
A formal distribution $a(z,w)\in U[[z,z^{-1},w,w^{-1}]]$ is local if and only if $a(z,w)$ can be written as
\begin{equation*}
a(z,w)=\sum_{j\in\Z^+}c^j(w)\frac{\partial^j_w\delta(z,w)}{j!},
\end{equation*}
for some $c^j(w)\in U[[w,w^{-1}]]$.
\end{prop}

In this paper, we adopt the following definition of Lie conformal algebras using $\lambda$-brackets as in Ref.$7$.
\begin{defi}\label{D1}\rm
A Lie conformal algebra is a $\C[\partial]$-module $A$ endowed with a $\lambda$-bracket $[a{}\, _\lambda \, b]$ which defines a linear map $A\otimes A\rightarrow A[\lambda]$, where $\lambda$ is an indeterminate and $A[\lambda]=\C[\lambda]\otimes A$, subject to the following axioms:
\begin{equation}\label{conformal}
\aligned
&[\partial a\,{}_\lambda \,b]=-\lambda[a\,{}_\lambda\, b],\ \ \ \
[a\,{}_\lambda \,\partial b]=(\partial+\lambda)[a\,{}_\lambda\, b];\\
&[a\, {}_\lambda\, b]=-[b\,{}_{-\lambda-\partial}\,a];\\
&[a\,{}_\lambda\,[b\,{}_\mu\, c]]=[[a\,{}_\lambda\, b]\,{}_{\lambda+\mu}\, c]+[b\,{}_\mu\,[a\,{}_\lambda \,c].
\endaligned
\end{equation}
\end{defi}

For any local formal distribution $a(z,w)$, the formal Fourier transformation $F^\lambda_{z,w}$ is defined as
$F^\lambda_{z,w}a(z,w)={\rm Res}_z e^{\lambda(z-w)}a(z,w)$.
Suppose $(g, F)$ is a formal distributions Lie algebra. Define $[a(z)\,{}_\lambda\, b(z)]=F^\lambda_{z,w}[a(z),b(w)]$ for any $a(z),b(z)\in F$. One can check  that this definition of $\lambda$-brackets satisfies (\ref{conformal}).
Given a formal distributions Lie algebra $(g, F)$, we may always include $F$ in the minimal family $F^c$ of pairwise local distributions which is closed under
the derivative $\partial$ and the $\lambda$-brackets.
Then $F^c$ is actually a Lie conformal algebra.
We call it a Lie conformal algebra of $g$.

Following this procedure,
we define $L_i(z)=\sum_{\a\in\Z} L_{\a,i}z^{-\a-2}$ for any $i\in\Z$.
One can show easily that
\begin{equation*}
[L_i(z),L_j(w)]=-\delta(z,w)\partial_w L_{i+j}(w)-2L_{i+j}(w)\partial_w\delta(z,w).
\end{equation*}
Hence by Proposition \ref{local}, $\{L_i(z)\,|\,i\in\Z\}$ are pairwise local formal distributions, whose coefficients span $\L$.
One can easily obtain that $[L_i\,{}_\lambda\, L_j]=-\partial L_{i+j}-2\lambda L_{i+j}$
through some direct computations. Under this definition of $\lambda$-brackets, $\mathscr{CW}:=\oplus_{i\in\Z} \C[\partial]{L_i}$ is a Lie conformal algebra of $\L$. We call $\mathscr{CW}$ {\it loop Virasoro Lie conformal algebra} in this paper.


\begin{defi}\label{D11}\rm
Let $A$ be a Lie conformal algebra. A linear map $\phi_\lambda: A\rightarrow A[\lambda]$ is called a conformal derivation if the following equalities hold:
\begin{equation*}
\aligned
&\phi_\lambda(\partial v)=(\partial+\lambda)\phi_\lambda(v),\ \ \ \
&\phi_\lambda([a\,{}_\mu \,b])=[(\phi_\lambda a)\,{}_{\lambda+\mu} \,b]+[a\,{}_\mu \,(\phi_\lambda b)].
\endaligned
\end{equation*}
We often write $\phi$ instead of $\phi_\lambda$ for simplicity.
\end{defi}
It can be easily verified that for any $x\in A$, the map ${\rm ad}_x$, defined by $({\rm ad}_x)_\lambda y= [x\, {}_\lambda\, y]$ for   $y\in A$, is a conformal derivation of $A$.
All conformal derivations of this kind are called {\it inner}.
Denote by ${\rm Der\,}A$
and ${\rm Inn\,}A$ the vector spaces of all conformal derivations and inner conformal derivations of $A$, respectively.
\begin{defi}\label{D2}\rm
A conformal module $M$ over a Lie conformal algebra $A$ is a $\C[\partial]$-module endowed with a $\lambda$-action $A\otimes M\rightarrow M[\lambda]$ such that
\begin{equation}
\aligned
&(\partial a)\,{}_\lambda\, v=-\lambda a\,{}_\lambda\, v,\ \ \ \ \ a{}\,{}_\lambda\, (\partial v)=(\partial+\lambda)a\,{}_\lambda\, v;\\
&a\,{}_\lambda\, (b{}\,_\mu\, v)-b\,{}_\mu\,(a\,{}_\lambda\, v)=[a\,{}_\lambda\, b]\,{}_{\lambda+\mu}\, v.
\endaligned
\end{equation}
\end{defi}
\begin{defi}\rm
A Lie conformal algebra $A$ is {\it $\Z$-graded} if $A=\oplus_{i\in \Z}A_i$ 
and $[A_i\,{}_\lambda\, A_j]\subset A_{i+j}[\lambda]$ for any $i,j\in \Z$.
Similarly, a conformal module $V$ over $A$ is  {\it $\Z$-graded} if $V=\oplus_{i\in \Z}V_i$ and $A_i\,{}_\lambda\, V_j\subset V_{i+j}[\lambda]$ for any $i,j\in \Z$. In addition, if each $V_i$ is freely generated by one element $v_i\in V_i$ over $\C[\partial]$, we call $V$ a {\it $\Z$-graded free intermediate series module}.
\end{defi}

\section{CONFORMAL DERIVATIONS OF $\mathscr{CW}$}

Assume $D\in {\rm Der} \mathscr{CW}$.
Define $D^i(L_j)=\pi_{i+j} D(L_j)$ for any $j\in\Z$, where in general $\pi_{i}$ is the natural projection from $\C[\lambda]\otimes \mathscr{CW}\cong \oplus_{k\in\Z}\C[\partial,\lambda]L_k$ onto $\C[\partial,\lambda]{L_{i}}$.
Then $D^i$ is a conformal derivation and $D=\sum_{i\in\Z} D^i$ in the sense that for any $x\in \mathscr{CW}$ only finitely many $D^i_\lambda(x)\neq0$.
\begin{lemm}\label{lc1}
For any $c\in\Z$, $D^c$ is an inner conformal derivation of the form $D={\rm ad}_{{g(\partial)}L_c}$ for some $g(\partial)\in\C[\partial]$.
\end{lemm}
\ni\ni{\it Proof.}\ \   Assume $D^c_\lambda(L_i)=f_i(\partial,\lambda)L_{i+c}$. Let $D^c_\lambda$ act on $[L_0\ {}_\mu \ L_i]=(-\partial-2\mu)L_i$,
one has
\begin{equation}\label{c1}
(-\partial-\lambda-2\mu)f_i(\partial,\lambda)=(-\partial-2\lambda-2\mu)f_0(-\lambda-\mu,\lambda)+(-\partial-2\mu)f_i(\partial+\mu,\lambda).
\end{equation}
Setting $\mu=0$ in (\ref{c1}), one gets
\begin{equation}
\lambda f_i(\partial,\lambda)=(\partial+2\lambda)f_0(-\lambda,\lambda).
\end{equation}
Hence $\lambda$ is a factor of  $f_0(-\lambda,\lambda)$ in the polynomial ring $C[\partial,\lambda]$. Setting $g(\lambda)=\frac{f_0(\lambda,-\lambda)}{\lambda}$,
we have $D={\rm ad}_{{g(\partial)}L_c}$. \QED
\begin{prop}
${\rm Der}\mathscr{CW}={\rm Inn}\mathscr{CW}$.
\end{prop}
\ni\ni{\it Proof.}\ \ From Lemma \ref{lc1}, we have $D=\sum_{i\in\Z} D^i=\sum_{i\in\Z} {\rm ad}_{{h_i(\partial)}L_i}$ for some $h_i(\partial)\in\C[\partial]$. If $h_i(\partial)\ne0$ for infinite many $i$'s, then $D(L_0)=\sum_{i\in\Z}D^i_\lambda(L_0)=\sum_{i\in\Z}[h_i(\partial)L_i\,{}_\lambda \, L_0]=\sum_{i\in\Z}(-\partial-2\lambda)h_i(-\lambda)L_i$ is an infinite sum, a contradiction to the definition of derivations. Thus $D=\sum_{i\in\Z} {\rm ad}_{{h_i(\partial)}L_i}={\rm ad}_h$ is a finite sum, where $h=\sum_{i\in\Z}h_i(\partial)L_i\in\mathscr{CW}$, i.e., $D\in {\rm Inn}\mathscr{CW}$.\QED

\section{RANK ONE CONFORMAL MODULES OVER $\mathscr{CW}$}

Suppose  $M$ is a free conformal module of rank one over $\mathscr{CW}$.
We may write $M=\C[\partial]v$ and assume $L_i\,{}_\lambda\, v=f_i(\partial,\lambda)v$, where $f_i(\partial,\lambda)\in\C[\partial,\lambda]$.
One can obtain the following result easily.
\begin{lemm}\label{lfo1}
For any $i,j\in\Z$, the following equality holds{\rm:}
\begin{equation}\label{m0}
(\mu-\lambda)f_{i+j}(\partial,\lambda+\mu)=f_j(\partial+\lambda,\mu)f_i(\partial,\lambda)-f_i(\partial+\mu,\lambda)f_j(\partial,\mu).
\end{equation}
\end{lemm}
\ni\ni{\it Proof.}\ \ A direct computation shows that
\begin{eqnarray*}&&
[L_i\,{}_\lambda\, L_j]\,{}_{\lambda+\mu} v=(-\partial L_{i+j}-2\lambda L_{i+j}){}_{\lambda+\mu} v=(\mu-\lambda)f_{i+j}(\partial,\lambda+\mu)v,
\\&& 
L_i\,{}_\lambda\,(L_j\, {}_\mu v)=L_i\,{}_\lambda (f_j(\partial,\mu)v)=f_j(\partial+\lambda,\mu)L_i {}_\lambda v=f_j(\partial+\lambda,\mu)f_i(\partial,\lambda)v,
\end{eqnarray*}
and
\begin{equation*}
L_j\,{}_\mu\, (L_i \,{}_\lambda v)=f_i(\partial+\mu,\lambda)f_j(\partial,\mu)v.
\end{equation*}
Then the result follows.\QED
\begin{lemm}\label{lfo2}
Suppose  $M$ is a nontrivial free conformal modules of rank one over $\mathscr{CW}$.
Then $f_0(\partial,\lambda)=a\lambda+b-\partial$ for some $a,b\in\C$.
\end{lemm}
\ni\ni{\it Proof.}\ \ Setting $j=0$ in (\ref{m0}), one has
\begin{equation}\label{m1}
(\mu-\lambda)f_{i}(\partial,\lambda+\mu)=f_0(\partial+\lambda,\mu)f_i(\partial,\lambda)-f_i(\partial+\mu,\lambda)f_0(\partial,\mu).
\end{equation}
Setting $i=0$ in (\ref{m1}), we have
\begin{equation}\label{m2}
(\mu-\lambda)f_{0}(\partial,\lambda+\mu)=f_0(\partial+\lambda,\mu)f_0(\partial,\lambda)-f_0(\partial+\mu,\lambda)f_0(\partial,\mu).
\end{equation}
Comparing the highest degree of $\lambda$ on both sides in (\ref{m2}),
\begin{equation*}
1+{\rm deg}_2 f_{0}={\rm deg}_1 f_{0}+{\rm deg}_2 f_{0},
\end{equation*}
where ${\rm deg}_1 f_{0}$ (or ${\rm deg}_2 f_{0}$) denotes the maximal degree of the first (or second) variable $x$ (or $y$) in the polynomial $f_{0}(x,y)$.
Thus ${\rm deg}_1 f_0=1$ and we can assume $f_0(\partial,\lambda)=c(\lambda)\partial+d(\lambda)$ for some $c(\lambda), d(\lambda)\in\C[\lambda]$.
Then we can rephrase (\ref{m2}) as follows:
\begin{equation}\label{m3}
\aligned
&(\mu-\lambda)(c(\lambda+\mu)\partial+d(\lambda+\mu))\\
=&(c(\mu)(\partial+\lambda)+d(\mu))(c(\lambda)\partial+d(\lambda))
-(c(\lambda)(\partial+\mu)+d(\lambda))(c(\mu)\partial+d(\mu)).
\endaligned
\end{equation}
Thus we obtain
\begin{equation}\label{m4}
c(\lambda+\mu)=-c(\lambda)c(\mu).
\end{equation}
and
\begin{equation}\label{m5}
(\mu-\lambda)d(\lambda+\mu)=c(\mu)d(\lambda)\lambda-c(\lambda)d(\mu)\mu.
\end{equation}
Then $c(\lambda)=0$ or $-1$ from (\ref{m4}).

If $c(\lambda)=0$, then $d(\lambda)=0$ from (\ref{m5}).
Consequently, $f_0(\partial,\lambda)=0$ and $f_i(\partial,\lambda)=0$ for any $i\in\Z$ from (\ref{m1}).
In this case, $M$ is a trivial $\mathscr{CW}$-module.

If $c(\lambda)=-1$, then $\mu d(\mu)-\lambda d(\lambda)=(\mu-\lambda)d(\lambda+\mu)$ from (\ref{m5}).
Consequently, $d(\lambda)=a\lambda+b$ for some $a, b\in\C$. Thus $f_0(\partial,\lambda)=a\lambda+b-\partial$.\QED

Next we attempt to compute $f_{i}(\partial,\lambda)$ for any $i\in\Z$.
Setting $\lambda=0$ in (\ref{m1}), we have
\begin{equation}\label{m6}
\mu f_{i}(\partial,\mu)=(f_i(\partial,0)-f_i(\partial+\mu,0))f_0(\partial,\mu).
\end{equation}
Since $f_0(\partial,\lambda)=-\partial+a\lambda+b$ is an irreducible polynomial in $\C[\partial,\lambda]$, which is a unique
factorization domain,
one obtain from \eqref{m6} that $(-\partial+a\lambda+b)\, | \,f_{i}(\partial,\lambda)$.
Thus we can write $f_{i}(\partial,\lambda)=(-\partial+a\lambda+b)g_i(\partial,\lambda)$ for some $g_i(\partial,\lambda)\in\C[\partial,\lambda]$.
Then we have
\begin{equation}\label{m7}
\mu g_{i}(\partial,\mu)=(-\partial+b)g_i(\partial,0)-(-\partial-\mu+b)g_i(\partial+\mu,0).
\end{equation}
Setting $\lambda=\mu$ in (\ref{m1}), we have
\begin{equation}\label{m8}
f_0(\partial+\lambda,\lambda)f_i(\partial,\lambda)=f_i(\partial+\lambda,\lambda)f_0(\partial,\lambda).
\end{equation}
Then we obtain
\begin{equation}\label{m9}
g_i(\partial+\lambda,\lambda)=g_i(\partial,\lambda)\mbox{ \ for \ }i\in\Z.
\end{equation}
Write $g_i(\partial,\lambda)=\sum_{j=0}^{n} a_j(\lambda)\partial^j$ for some $a_j(\lambda)\in\C[\partial]$.
From (\ref{m9}), we obtain
\begin{equation}\label{m10}
g_i(k\lambda,\lambda)=g_i(0,\lambda)=a_0(\lambda)\mbox{ \ for \ }k\in\N.
\end{equation}
Take $k$ to be integers from $1$ to $n$, we get a series of linear equations in unknowns $a_j(\lambda)\lambda^j$ for $j=1,\cdots,n$,
whose coefficients form a nondegenerate Vandermonde matrix.
So $a_j(\lambda)\lambda^j=0$ for $j=1,\cdots,n$.
Now we have shown that $g_i(\partial,\lambda)$ is actually a polynomial in the second variable.
Then from (\ref{m7}), $g_i(\lambda)$ is actually a constant.
Hence we can assume that $f_i(\partial,\lambda)=c_i f_0 (\partial,\lambda)$ for some $c_i\in\C$.
Under this assumption, we have $c_{i+j}=c_i c_j$ for  $i, j\in\Z$.
Hence $c_i=c_1^i$ for $i\in\Z$.

In this case, the conformal module structure on $M$ is described by three parameters $a,b,c\in\C$, where the $\lambda$-action is given by
\begin{equation}\label{act000}
L_i\, {}_\lambda \,v=c^i(-\partial+a\lambda+b)v\mbox{ \  for \ }i\in\Z.\end{equation}
Denote this module by $V_{a,b,c}$.

Now we get the main result of this section.
\begin{prop}
Assume $V$ is a nontrivial free module of rank one. Then $V$ is isomorphic to $V_{a,b,c}$ defined by \eqref{act000} for some $a,b,c\in\C$.\QED
\end{prop}

\section{$\Z$-GRADED FREE INTERMEDIATE SERIES MODULES OVER $\mathscr{CW}$}

In this section, we give a classification of $\Z$-graded free intermediate series modules over $\mathscr{CW}$.
Let $V$ be an arbitrary $\Z$-graded free intermediate series module over $\mathscr{CW}$. Then $V=\oplus_{i\in\Z}V_i$, where each $V_i$
is freely generated by some element $v_i\in V_i$ over $\C[\partial]$.
For any $i,j\in\Z$, denote $L_i\,{}_\lambda\, v_j=f_{i,j}(\partial,\lambda)v_{i+j}$.
We call $\{f_{i,j}(\partial,\lambda)\,|\,i,j\in\Z\}$ the {\it structure
coefficients} of $V$ related to the $\C[\partial]$-basis $\{v_i\,|\,i\in\Z\}$. Then the conformal module structure on $V$ is determined
if and only if all of its structure coefficients are specified.

From Definition \ref{D2}, one can obtain the following result.
\begin{lemm}\label{lf1}
The structure coefficients $f_{i,j}(\partial,\lambda)$ of
$V$ satisfy the following equality {\rm :}
\begin{equation}\label{f1}
(\mu-\lambda)f_{i+j,k}(\partial,\lambda+\mu)=f_{j,k}(\partial+\lambda,\mu)f_{i,j+k}(\partial,\lambda)-f_{i,k}(\partial+\mu,\lambda)f_{j,i+k}(\partial,\mu).
\end{equation}
\end{lemm}
\ni\ni{\it Proof.}\ \   One can check it easily that
\begin{equation*}
[L_i\,{}_\lambda\, L_j]\,{}_{\lambda+\mu} \,v_k=(-\partial L_{i+j}-2\lambda L_{i+j})\,{}_{\lambda+\mu} \,v_k=(\mu-\lambda)f_{i+j,k}(\partial,\lambda+\mu)v_{i+j+k}.
\end{equation*}
Similarly we  have the following two equalities:
\begin{equation*}
\aligned
L_i\,{}_\lambda \,(L_j\, {}_\mu \,v_k)&=L_i\,{}_\lambda \,(f_{j,k}(\partial,\mu)v_{j+k})
=f_{j,k}(\partial+\lambda,\mu)f_{i,j+k}(\partial,\lambda)v_{i+j+k};\\
L_j\,{}_\mu \,(L_i \,{}_\lambda\, v_k)&=f_{i,k}(\partial+\mu,\lambda)f_{j,i+k}(\partial,\mu)v_{i+j+k}.
\endaligned
\end{equation*}
Then the result follows from the defining relations of a conformal module.\QED

\begin{lemm}\label{lf2}
If $f_{j,k}(\partial,\lambda)=0$ for some $j,k\in\Z$, then $f_{i+j,k}(\partial,\lambda)=0$ for all $i\in\Z$.
\end{lemm}
\ni\ni{\it Proof.}\ \  From (\ref{f1}), we get
\begin{equation}
(\mu-\lambda)f_{i+j,k}(\partial,\lambda+\mu)=-f_{i,k}(\partial+\mu,\lambda)f_{j,i+k}(\partial,\mu).
\end{equation}
Then $(\mu-\lambda)\,|\,f_{i,k}(\partial+\mu,\lambda)f_{j,i+k}(\partial,\mu)$. Since $(\mu-\lambda)$ is an irreducible polynomial in $\C[\partial,\lambda,\mu]$, which is a unique factorization domain, $(\mu-\lambda)\,|\,f_{i,k}(\partial+\mu,\lambda)$. Then $f_{i,k}(\partial+\lambda,\lambda)=0$. Thus $f_{i,k}(\partial,\lambda)=0$. Consequently, $f_{i+j,k}(\partial,\lambda+\mu)=0$. \QED

\begin{lemm}\label{lf3}
If the structure coefficient $f_{0,j+k}(\partial,\lambda)=0$, then $f_{j,k}(\partial,\lambda)=0$.
\end{lemm}
\ni\ni{\it Proof.}\ \  Let $i=0$ in (\ref{f1}), then we obtain
\begin{equation}\label{f7}
(\mu-\lambda)f_{j,k}(\partial,\lambda+\mu)=f_{j,k}(\partial+\lambda,\mu)f_{0,j+k}(\partial,\lambda)-f_{0,k}(\partial+\mu,\lambda)f_{j,k}(\partial,\mu).
\end{equation}
If $f_{0,j+k}(\partial,\lambda)=0$, then
\begin{equation*}
(\mu-\lambda)f_{j,k}(\partial,\lambda+\mu)=-f_{0,k}(\partial+\mu,\lambda)f_{j,k}(\partial,\mu).
\end{equation*}
Thus $f_{j,k}(\partial,\lambda)=0$.\QED

\begin{prop}\label{pf1}
If $f_{j_0,k_0}(\partial,\lambda)=0$ for some $j_0,k_0\in\Z$, then $f_{j,k}(\partial,\lambda)=0$ for all $j,k\in\Z$.
\end{prop}
\ni\ni{\it Proof.}\ \ It follows from Lemmas \ref{lf2} and \ref{lf3} immediately.\QED

From now on, we assume that $V$ is a nontrivial module. Then from Proposition \ref{pf1}, $f_{j,k}(\partial,\lambda)\neq0$ for all $j,k\in\Z$.

\begin{lemm}\label{lf4}
For any $k\in\Z$, there exist $a_k,b_k\in\C$ such that
the structure coefficient $f_{0,k}(\partial,\lambda)$ is of the form $a_k \lambda+b_k-\partial$.
\end{lemm}
\ni\ni{\it Proof.}\ \   Set $i=j=0$ in (\ref{f1}), one gets
\begin{equation}\label{f2}
(\mu-\lambda)f_{0,k}(\partial,\lambda+\mu)=f_{0,k}(\partial+\lambda,\mu)f_{0,k}(\partial,\lambda)-f_{0,k}(\partial+\mu,\lambda)f_{0,k}(\partial,\mu).
\end{equation}
Then the statement follows from an argument analogous to that of Lemma \ref{lfo2}.\QED

\begin{lemm}\label{lf5}
For any $j,k\in\Z$, $b_k=b_{j+k}$.
\end{lemm}
\ni\ni{\it Proof.}\ \ Setting $\lambda=0$ in (\ref{f7}), one has
\begin{equation}\label{f8}
\mu f_{j,k}(\partial,\mu)=f_{j,k}(\partial,\mu)(b_{j+k}-\partial)-(b_k-\partial-\mu)f_{j,k}(\partial,\mu).
\end{equation}
If $b_k\neq b_{j+k}$, we have $f_{j,k}(\partial,\mu)=0$. A contradiction.\QED

Write $b=b_k$, then we have $f_{0,k}(\partial,\lambda)=a_k \lambda+b-\partial$ by Lemma \ref{lf4}.
Now we can prove the following result, which plays a key role in the classification of $\Z$-graded free intermediate series modules.
\begin{lemm}\label{lf6}
For any $j,k\in\Z$, there exists $c_{j,k}\in\C^*$ such that one and only one of the following three equalities holds:
\begin{equation*}
\aligned
&(1)\ \ f_{j,k}(\partial,\lambda)=c_{j,k}(\partial-b-a_k\lambda),\ \ \ a_k=a_{j+k};\\
&(2)\ \ f_{j,k}(\partial,\lambda)=c_{j,k},\ \ \ (a_k,a_{j+k})=(0,-1);\\
&(3)\ \ f_{j,k}(\partial,\lambda)=c_{j,k}(\partial-b)(\partial-b+\lambda),\ \ \ (a_k,a_{j+k})=(-1,0).\\
\endaligned
\end{equation*}
\end{lemm}
\ni\ni{\it Proof.}\ \ Let $\mu=0$ in (\ref{f7}), we have
\begin{equation}\label{f9}
-\lambda f_{j,k}(\partial,\lambda)=f_{j,k}(\partial+\lambda,0)(a_{j+k} \lambda+b-\partial)-(a_k \lambda+b-\partial-\mu)f_{j,k}(\partial,0).
\end{equation}
Since $\lambda$ is a factor of $f_{j,k}(\partial+\lambda,0)-f_{j,k}(\partial,0)$ in $\C[\partial,\lambda]$,
we can rephrase the equation (\ref{f9}) as follows
\begin{equation}\label{f10}
f_{j,k}(\partial,\lambda)=(\partial-b)\frac{f_{j,k}(\partial+\lambda,0)-f_{j,k}(\partial,0)}{\lambda}+a_k f_{j,k}(\partial,0)-a_{j+k}f_{j,k}(\partial+\lambda,0).
\end{equation}
Denote $d(\partial)=f_{j,k}(\partial,0)$ and let $\lambda$ approach to zero in (\ref{f10}), then we obtain an ordinary differential equation (ODE):
\begin{equation*}
(1+a_{j+k}-a_k)d(\partial)=(\partial-b)d'(\partial).
\end{equation*}
As $d(\partial)$ is a polynomial, there exists some $c_{j,k}\in\C$ such that
\begin{equation}\label{f12}
d(\partial)=c_{j,k}(\partial-b)^{(1+a_{j+k}-a_k)}.
\end{equation}
Since $f_{j,k}(\partial,\lambda)\neq0$, then $c_{j,k}\neq0$ and $1+a_{j+k}-a_k\in\N$.
Combining (\ref{f10}) and (\ref{f12}), one gets
\begin{equation}\label{f13}
\aligned
f_{j,k}(\partial,\lambda)=&c_{j,k}a_k(\partial-b)^{(1+a_{j+k}-a_k)}-c_{j,k}a_{j+k}(\partial+\lambda-b)^{(1+a_{j+k}-a_k)}\\
&+c_{j,k}(\partial-b)\lambda^{-1}[(\partial+\lambda-b)^{(1+a_{j+k}-a_k)}-(\partial-b)^{(1+a_{j+k}-a_k)}].
\endaligned
\end{equation}
Let $\lambda=\mu$ in (\ref{f7}), we have
\begin{equation}\label{f14}
f_{j,k}(\partial+\lambda,\lambda)(a_{j+k} \lambda+b-\partial)=(a_k \lambda+b-\partial-\lambda)f_{j,k}(\partial,\lambda).
\end{equation}
Let $\partial-b=N\lambda$ for some $N\in\C$ in (\ref{f14}), one can obtain by (\ref{f13}) that
\begin{equation}\label{f15}
\aligned
&(N+2)^e(N+1-a_{j+k})(a_{j+k}-N)\\
=&(a_k-1-N)(a_{k}-N)N^e+2(a_k-1-N)(N-a_{j+k})(N+1)^e,
\endaligned
\end{equation}
where $e=1+a_{j+k}-a_k$.
Let $N=-1$ in (\ref{f15}), we have
\begin{equation*}
-a_{j+k}(a_{j+k}+1)=a_k (a_{k}+1) {(-1)}^{(1+a_{j+k}-a_k)}.
\end{equation*}
If $1+a_{j+k}-a_k$ is odd, then
$a_{j+k}(a_{j+k}+1)=a_k (a_{k}+1)$.
Thus $a_k=a_{j+k}$ or $a_k+a_{j+k}+1=0$.
If $a_k=a_{j+k}$, then $f_{j,k}(\partial,\lambda)=-c_{j,k}(a_k \lambda+b-\partial)$.
If $a_k+a_{j+k}+1=0$, then $a_k=-\frac{2n+1}{2}$ and $a_{j+k}=\frac{2n-1}{2}$ for some $n\in\N$.
Set $N=0$ in (\ref{f15}), we have
\begin{equation*}
2^e(1-a_{j+k})a_{j+k}=-2(a_k-1)a_{j+k}.
\end{equation*}
In the case of $a_k=-\frac{2n+1}{2}$ and $a_{j+k}=\frac{2n-1}{2}$, we have $2^{2n}(3-2n)=2n+3$.
Since $n\in\N$, then $n=0$. Consequently $a_k=a_{j+k}=-\frac{1}{2}$ and return to the first case.

If $1+a_{j+k}-a_k$ is even, then
$a_{j+k}(a_{j+k}+1)=-a_k (a_{k}+1)$.
Thus $(a_{j+k}+\frac{1}{2})^2+(a_{k}+\frac{1}{2})^2=\frac{1}{2}$.
So $(a_k,a_{j+k})=(0,-1)$ or $(a_k,a_{j+k})=(-1,0)$.
In the case of $(a_k,a_{j+k})=(0,-1)$, $f_{j,k}(\partial,\lambda)=c_{j,k}$
In the case of $(a_k,a_{j+k})=(-1,0)$, $f_{j,k}(\partial,\lambda)=c_{j,k}(\partial-b)(\partial-b+\lambda)$.
It completes the proof. \QED

Let $A$ be the sequence $\{a_i\}_{i\in\Z}$. Then by Lemma \ref{lf6}, $a_i\in\{0,-1\}$ or $a_i=a_j$ for any $i,j\in\Z$.

We can prove the following result.
\begin{lemm}
For any $i,j,k\in\Z$, $-c_{i+j,k}=c_{j,k}c_{i,j+k}$.
\end{lemm}
\ni\ni{\it Proof.}\ \
If there exists a constant $a\in\C$ such that $a_i=a$ for all $i\in\Z$, then $f_{j,k}(\partial,\lambda)=-c_{j,k}(a \lambda+b-\partial)$.
In this case, we can rephrase (\ref{f1}) as follows:
\begin{equation*}
\aligned
&c_{i+j,k}(\lambda-\mu)(a \lambda+a \mu+b-\partial)\\
=&c_{j,k}c_{i,j+k}(a \mu+b-\partial-\lambda)(a \lambda+b-\partial)-c_{i,k}c_{j,i+k}(a \lambda+b-\partial-\mu)(a \mu+b-\partial).
\endaligned
\end{equation*}
Then $-c_{i+j,k}=c_{j,k}c_{i,j+k}$.

If for any $i\in\Z$, $a_i\in\{0,-1\}$. Then we can check the desired equality case by case from Lemma \ref{lf1}. \QED

Since $-c_{i+j,k}=(-c_{j,k})(-c_{i,j+k})$, we get the following result.$^4$
\begin{lemm}
There exist some nonzero numbers $d_k\in\C$ such that $-c_{j,k}=\frac{d_{j+k}}{d_k}$. \QED
\end{lemm}
Let $u_j=d_j v_j$ for any $j\in\Z$. Then $\{u_j\,|\,j\in\Z\}$ is also a $\C[\partial]-$basis of $V$. Replacing $\{v_j\,|\,j\in\Z\}$ by $\{u_j\,|\,j\in\Z\}$, we can assume $c_{j,k}=-1$ for any $j,k\in\Z$.

Given $a,b\in\C$, we can construct a nontrivial $\Z$-graded free intermediate series module $V_{a,b}$, where $V_{a,b}=\oplus_{i\in\Z}\C[\partial]v_i$ and
the $\lambda$-actions are given by
\begin{equation}\label{act1111}
L_i\,{}_\lambda\, v_j=-(\partial-b-a\lambda)v_{i+j}.
\end{equation}

Let $A=\{a_i\}_{i\in\Z}$ be a sequence with $a_i\in\{0,-1\}$ for any $i\in\Z$ and $b\in\C$.
We can construct a nontrivial $\Z$-graded free intermediate series module $V_{A,b}$, where $V_{A,b}=\oplus_{i\in\Z}\C[\partial]v_i$ and
the $\lambda$-actions are given by
\begin{equation}\label{act11112}
L_i\,{}_{{}_\lambda} v_j=
\begin{cases}
-(\partial-b)v_{i+j}, &\ if\  (a_j,a_{i+j})=(0,0);\\[4pt]
-(\partial-b+\lambda)v_{i+j},&\  if \ (a_j,a_{i+j})=(-1,-1);\\[4pt]
-v_{i+j}, &\ if \ (a_j,a_{i+j})=(0,-1);\\[4pt]
-(\partial-b)(\partial-b+\lambda)v_{i+j},&\  if \ (a_j,a_{i+j})=(-1,0).
\end{cases}
\end{equation}

Then we get the main result of this section.
\begin{theo}
Assume $V$ is a nontrivial $\Z$-graded free intermediate series module over $\mathscr{CW}$.
Then $V$ is isomorphic to $V_{a,b}$ defined by \eqref{act1111} or $V_{A,b}$ defined by \eqref{act11112}.\QED
\end{theo}

\ni\ni\textbf{ACKNOWLEDGMENTS}

This work was supported by NSF Grant No.11371278, 11001200, 11271284 and 11101269 of China, the Grant No.12XD1405000 of Shanghai Municipal Science and Technology Commission and the Fundamental Research Funds for the Central Universities.

\end{document}